\documentclass[11pt]{article}\usepackage[mathletters]{ucs}\usepackage[utf8x]{inputenc}%

\usepackage{color}
\usepackage{graphicx}
\usepackage{amsmath}
\usepackage{amsthm}
\usepackage{amssymb}
\usepackage{amsxtra}
\usepackage{amsfonts}

\theoremstyle{plain}
\newtheorem{theorem}{Theorem}

\theoremstyle{definition}

\theoremstyle{remark}

\newcommand{\myref}[1]{\ref{#1}}
\newcommand{\myeqref}[1]{\eqref{#1}}
\newcommand{\myeqlab}[1]{\label{#1}}

\def\[#1\]{\begin{alignat}{5}#1\end{alignat}}

\def\eqdef{=:}
\newcommand{\half}{\frac12}
\newcommand{\const}{\text{const}}
\newcommand{\defm}[1]{\emph{#1}}
\newcommand{\conv}{\rightarrow}
\newcommand{\avg}[1]{\langle#1\rangle}
\newcommand{\jmp}[1]{[#1]}

\newcommand{\csep}{\quad,\quad}

\newcommand{\boi}[2]{{]#1,#2[}}

\newcommand{\subeq}[2]{\mathord{\underbrace{\mathop{#1}}_{#2}}}

\newcommand{\topref}[2]{\overset{\text{\myeqref{#1}}}{#2}}
\newcommand{\impl}{\Rightarrow}
\newcommand{\eqv}{\Leftrightarrow}
\newcommand{\sign}{\operatorname{sign}}

\newcommand{\crossp}{\times}
\newcommand{\dotp}{\cdot}

\newcommand{\xsv}{x}
\renewcommand{\vec}[1]{\mathbf{#1}}
\newcommand{\vv}{\vec\vsym}

\newcommand{\gisen}{γ}
\newcommand{\vsym}{u}
\newcommand{\vx}{\vsym^x}
\newcommand{\vy}{\vsym^y}

\newcommand{\vxu}{\vx_0}
\newcommand{\vyu}{\vy_0}
\newcommand{\csnd}{c}
\newcommand{\csndu}{\csnd_0}
\newcommand{\vnq}{(\vsym^n)^2}

\newcommand{\vvu}{\vv_0}

\newcommand{\vt}{\vsym^t}

\newcommand{\nn}{\vec n}
\newcommand{\ts}{\vec t}

\newcommand{\fod}{f}
\newcommand{\xiv}{\xi}
\newcommand{\etav}{\eta}
\newcommand{\xinor}{\xi_n}
\newcommand{\ximax}{\xi_M}

\newcommand{\Machu}{M_0}
\newcommand{\pp}{P}
\newcommand{\ppu}{P_0}
\newcommand{\dens}{\varrho}
\newcommand{\densu}{\dens_0}

\newcommand{\idens}{V}
\newcommand{\idensu}{\idens_0}
\newcommand{\epm}{e}
\newcommand{\epmu}{\epm_0}
\newcommand{\Hpm}{H}

\newcommand{\Epm}{E}

\newcommand{\hpm}{h}

\newcommand{\hpmu}{\hpm_0}
\newcommand{\spm}{S}

\newcommand{\spmu}{\spm_0}

\newcommand{\temp}{T}
\newcommand{\tempu}{T_0}

\newcommand{\jsym}{j}
\newcommand{\jj}{\vec\jsym}

\newcommand{\jn}{\jsym^n}
\newcommand{\jnq}{(\jsym^n)^2}
\newcommand{\jnu}{\jsym^n_0}

\newcommand{\Rspec}{R}
\newcommand{\cvspec}{c_v}
\newcommand{\cpspec}{c_p}
\newcommand{\turn}{\theta}
\newcommand{\vort}{\omega}
\newcommand{\vpot}{\phi}

\begin{document}

\title{Convexity of shock polars}
\author{Volker Elling}
\date{}
\maketitle
\begin{abstract}
  We show that the shock polars of compressible full potential flow are strictly convex if the enthalpy per mass is a convex function of volume per mass, in particular when the sound speed is a nondecreasing function of density. Counterexamples are given for some cases that violate the enthalpy condition. For the full Euler equations with convex equation of state satisfying the ideal gas law, polars are strictly convex if heat capacity is constant, but counterexamples are given in variable cases, showing no useful generalizations are possible. 
\end{abstract}

\section{Introduction}

\newif\ifGPblacktext%
\GPblacktextfalse%

\begin{figure}
  \noindent\parbox{.56\textwidth}{%
    {\tiny\input{intropolar.xxx}}
    \caption{Full Euler shock polar (gray) for upstream Mach number $\Machu=1.3$ and constant ratio of heats $\gisen=7/5$.}
    \label{fig:intropolar}
  }\hfil\parbox{.4\textwidth}{%
    {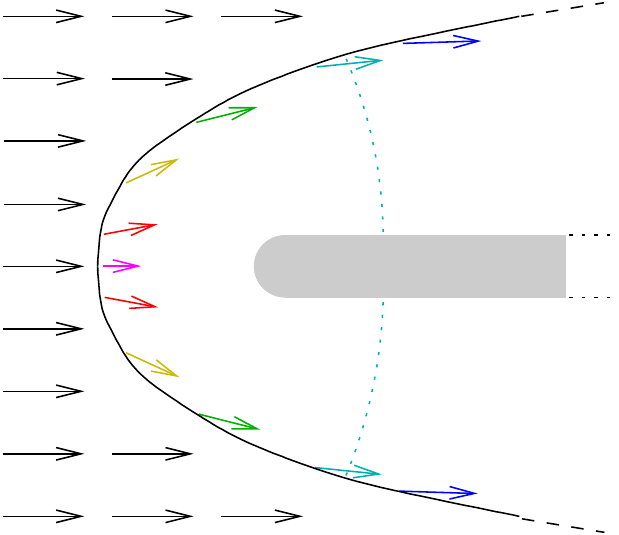}
    \caption{Bow shock in supersonic flow onto a blunt solid}
    \label{fig:bowshock}
  }%
\end{figure}

Consider supersonic flow onto a blunt solid body (fig.\ \myref{fig:bowshock}). Ahead of the body a detached bow shock forms, with constant velocity $\vvu$ on the upstream side. The downstream-side velocities $\vv$ vary with angle of shock to $\vvu$. The set of possible $\vv$ is called \defm{shock polar} (fig.\ \myref{fig:intropolar}). 
On the symmetry axis the shock is perpendicular to the inflow $\vvu$ (\defm{normal shock}). At \defm{critical shocks} the maximum possible angle between $\vvu$ and $\vv$ is realized. In the downstream infinity limit the shocks vanish, $\vv\conv\vvu$. 

\begin{figure}
  \noindent\hfil\parbox{.48\textwidth}{%
    \centerline{\footnotesize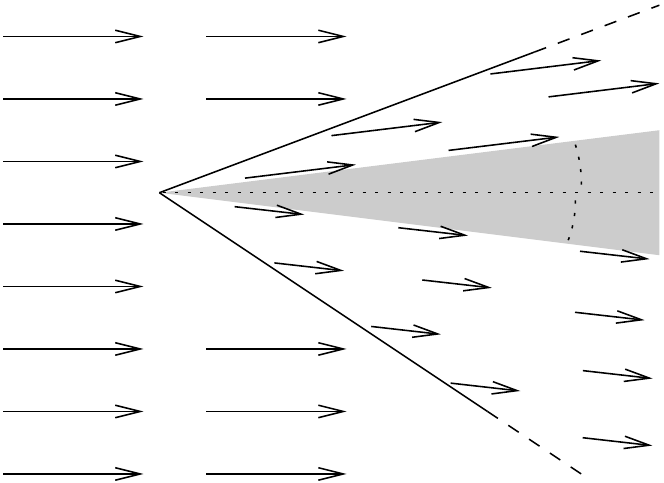}
    \caption{Supersonic flow onto a symmetric wedge producing a strong and a weak reflected shock (usually latter observed on both sides)}
    \label{fig:weakstrong}
  }\hfil\parbox{.48\textwidth}{%
    \centerline{\footnotesize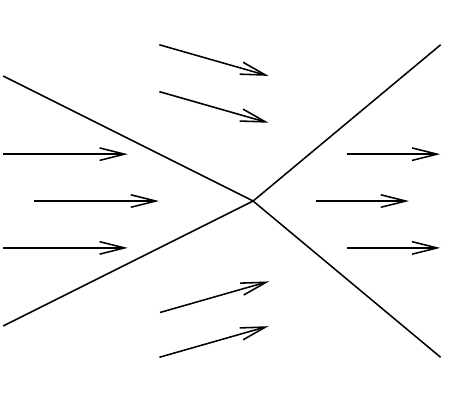}
    \caption{Two incident shocks producing two reflected shocks that turn velocity $\vv$ back to horizontal}
    \label{fig:fourshocks}
  }\hfil%
\end{figure}

Shock polars are also essential for studying shocks attached to solids, for example supersonic flow onto wedges or ramps (fig.\ \ref{fig:weakstrong}; \cite{yongqian-zhang-wedge-glimm,elling-liu-pmeyer}). The shock must turn the velocity by an angle $\turn$ keeping the flow parallel to the solid (\defm{slip condition}). For angles above critical no attached solutions are possible, for smaller angles there are usually two, called \defm{weak} and \defm{strong} (fig. \myref{fig:intropolar}). Interaction of shocks in interior flow in absence of any solids is also governed by shock polars; fig.\ \ref{fig:fourshocks} shows a simple case of two incoming and two outgoing shocks; for a large variety of more complicated Mach and Guderley reflections, multi-dimensional Riemann problems and other interactions can be found in \cite{yuxi-zheng-book,li-zhang-yang,geng-lai-indiana2019,shuxing-chen-book-2020,hunter-tesdall,eunheuikim-chungminlee-2013}. 

For the case of the full Euler equations with constant specific heat capacity, it is known classically (dating to the early 20th century) that each half of the polar has a unique critical shock, and that most weak shocks have supersonic $\vv$, whereas strong, critical and a small range of weak ones have subsonic $\vv$. Steady full Euler shocks satisfy in each point the following relations for conservation of mass, momentum and energy:
\[ 0 &= \jmp{\vv\dotp\nn\dens} ,
\\ 0 &= \jmp{\vv\dotp\nn \dens\vv + \pp\nn},
\\ 0 &= \jmp{\vv\dotp\nn\dens\Epm+\pp\vv\dotp\nn}. \]
Here $\jmp{f}=f-f_0$ is the difference between a downstream-side limit $f$ of some quantity and the upstream-side limit $f_0$ in the same point of the shock;
$\dens>0$ is mass density, $\vv$ the velocity vector field, $\nn$ a unit shock normal that points downstream, meaning mass flux $\jj=\dens\vv$ is positive in the $\nn$ direction, i.e.\ 
\[ \jnu = \densu\vvu\dotp\nn > 0 . \]
This distinguishes proper shocks from contact discontinuities across which there is no mass flux, $\jnu=0$.
$\Epm$ is total energy per mass, $\Epm=\epm+\half|\vv|^2$, with $\epm$ energy per mass. $\epm$ is linked to volume per mass $\idens=1/\dens$ and entropy per mass $\spm$ by some equation of state (eos) $\epm=\epm(\idens,\spm)$. 
The eos defines pressure
\[ \pp = - (\frac{\partial\epm}{\partial\idens})_\spm . \]
Here $(\partial X/\partial Y)_Z$ is the derivative of $X=X(Y,Z)$ with respect to $Y$ with second variable $Z$ held fixed. We omit $()_Z$ whenever clear from the context, in particular when $X=X(Y)$ has been assumed or shown to depend only on $Y$. For $\epm$ and its derivatives, $\partial_\idens$ is with $\spm$ held fixed, $\partial_\spm$ with $\idens$ held fixed. 
Temperature is defined by 
\[ \temp = (\frac{\partial\epm}{\partial\spm})_{\idens} , \]
whereas
\[ \cvspec = (\frac{\partial\epm}{\partial\temp})_{\idens} \]
is \defm{specific heat capacity at constant volume}, while
\[ \cpspec = (\frac{\partial\hpm}{\partial\temp})_{\pp} \]
is \defm{specific heat capacity at constant pressure},  
\[ \csnd = \sqrt{ (\frac{\partial\pp}{\partial\dens})_{\spm} } \]
is \defm{sound speed}, 
\[ \Hpm=\hpm+\frac12|\vv|^2 \myeqlab{eq:Hpm}\]
\defm{total enthalpy per mass}, for enthalpy per mass
\[ \hpm = \epm+\pp\idens . \]

Over wide ranges of temperature and pressure, many gases follow closely the \defm{ideal gas equation}
\[ \frac{\pp\idens}{\temp} = \Rspec = \const, \myeqlab{eq:idealgas}\]
where $\Rspec$ is the \defm{specific gas constant}. Those fluids are also called \defm{thermally perfect}.
It follows (as shown later) that for thermally perfect fluid $\epm=\epm(\temp)$ is a function of temperature alone (Joule's second law).

This still permits variable $\cvspec$, and indeed air and many other gases show significant variations with temperature, usually increase (see the introduction of \cite{elling-idealpolar} for a detailed discussion). In narrow temperature ranges it is reasonable to assume constant heat capacity $\cvspec$. In this case the ratio of heats (also called adiabatic exponent) $\gisen=\cpspec/\cvspec$ is often used, putting the eos in the form
\[ \epm = \epmu + \frac1{\gisen-1} (\frac{\exp\frac{\spm-\spmu}{\Rspec}}{\idens/\idensu})^{\gisen-1} \]
and
\[ \frac{\pp}{\ppu} = (\frac{\dens}{\densu})^{\gisen} \exp\big((\gisen-1)\frac{\spm-\spmu}{\Rspec}\big) . \]
Such fluids are called \defm{polytropic} (or \defm{calorically perfect}). 

After some calculation the full Euler shock relations take the convenient form
\[ 0 &= \jmp{\jn} ,                       \myeqlab{eq:full-jn}
\\ 0 &= \jmp{\vv\dotp\ts} ,               \myeqlab{eq:full-vt}
\\ 0 &= \jmp{\pp} + \jnq \jmp{\idens},    \myeqlab{eq:full-momn}
\\ 0 &= \jmp{\Hpm} .                      \myeqlab{eq:full-Hpm} \]
$\ts$ runs through all shock unit tangents; in two dimensions only one is needed which is chosen $90$ degree counterclockwise from the downstream unit normal $\nn$. 

If smooth full Euler solutions have zero vorticity $\vort=\nabla\crossp\vv$ at initial time, then also at later times. This leads to \defm{compressible potential flow}, which was studied in the early 20th century mostly in smooth regions. The full Euler equations in two dimensions are a system of four equations with the same number of scalar fields, whereas potential flow only uses a single scalar \defm{velocity potential} field $\vpot$ defined by $\vv=\nabla\vpot$. For this reason, in the second half of the century the need for transonic flow calculations on slow computers with small memory motivated the study of potential flow with shocks, either for the transonic small disturbance equations \cite{karman-tsd-1947,murman-cole,jameson-1978}, which will not be considered here, or for ``full'' potential flow, which uses a constant $\spm$ giving a scalar equation of state $\epm=\epm(\dens)$. w
Full potential flow shocks are defined by three of the four full Euler relations above:
\[ 0 &= \jmp{\jn}, \myeqlab{eq:potf-jn}
\\ 0 &= \jmp{\vv\dotp\ts}, \myeqlab{eq:potf-vt}
\\ 0 &= \jmp{\Hpm}. \myeqlab{eq:potf-Hpm}\] 
This is not to be confused with the \defm{barotropic Euler} equations (also called \defm{isentropic} or \defm{adiabatic}), which also use constant $\spm$, but keep \myeqref{eq:full-momn} instead of \myeqref{eq:full-Hpm}; that choice does not preserve zero vorticity after fluid passes through shocks.

As computers became cheaper and faster, interest on the numerical side decreased in favor of calculations with full Euler and more complicated models. On the other hand, a general trend on the PDE/analysis side towards multi-dimensional flow problems increased interest in potential flow there, since it allows approaching the complications of shock reflection while neglecting, in a first approach, the additional complications from possibly singular vorticity transport, which are not particular to compressible flow. As shock strength $ε=\dens/\densu-1$ decreases to $0$, entropy jump $\spm-\spmu$ decreases like $O(ε^3)$, so that potential flow shocks of small strength are very close to full Euler shocks. 

Compressible full potential flow also results from certain reformulations and limits of the nonlinear Schr\"odinger equation. In that case the potential $V(|\psi|)$ determines the eos $\pp(\dens)$. Limits of some water wave models also produce potential flow. Eos that are ``unnatural'' for gas dynamics may be natural in other applications.

The purpose of this note is to clarify exactly when shock polars (fig.\ \ref{fig:intropolar}) are strictly convex. The motivation is twofold. First, some earlier work \cite{elling-polarpotf,elling-idealpolar} received criticism for omitting discussion of polar convexity. Second, convexity is used, along with other polar properties, in many recent articles, sometimes in rigorous proof steps, without references or in-place justification. 
It is frequently stated, or used implicitly, that critical shocks are unique (in each half of the polar), that there are a unique weak and a unique strong shock for sufficiently small turning angle $\turn$, that the strong and critical shocks are always subsonic, and that the shock polar is a strictly convex curve. 
Sometimes \cite{courant-friedrichs} is cited (see Chapter 121), but this textbook collecting classical work only covers the $\gisen>1$ case of full Euler shocks, not potential flow shocks, which are \emph{not} a special case of full Euler, unlike the smooth flows.
Likewise, \cite{elling-liu-pmeyer} does not prove, claim, or use convexity.
To our knowledge, convexity or uniqueness have not been discussed for potential flow prior to \cite[Theorem 1]{elling-sonic-potf}, which establishes the result in the important but special case $1<\gisen<\infty$. 
Besides, for any model it is not safe to assume any of these properties continue to hold when moving from $\gisen>1$ polytropic eos to general ones. 

In this note we generalize the known results to larger classes of eos, while providing counterexamples that obstruct all apparent natural generalizations.
For potential flow,
if we assume, in addition to standard conditions (below),
that $\hpm$ is a convex function of $\idens$, then we can show the compressive part of the $\vv$ shock polar is strictly convex.
In particular, the assumption is satisfied by any eos with \defm{monotone sound speed}, meaning $\csnd$ is a nondecreasing function of $\dens$.
In the polytropic case the precondition allows any $\gisen>0$. But the result does not extend to $\gisen<0$; for some such $\gisen$ and some non-polytropic truncations of such eos we plot numerical \emph{counterexamples}. This is remarkable because \cite{elling-polarpotf} proved uniqueness of critical shocks using only standard assumptions, which permit all $\gisen>-1$ in the polytropic case. This shows clearly that shock polar convexity is a more ``fragile'' property. 

For the full Euler case, we recall the classical explicit shock polar formula in the polytropic case, from which strict convexity of the polar is clear; since the literature usually omits the detailed calculation, it is given below for the sake of completeness. We give numerical counterexamples showing polar convexity need not hold for convex eos that are ideal but non-polytropic. Counterexamples are possible even when the sound speed $\csnd$ is a nondecreasing function of temperature $\temp$. Hence, for full Euler, convexity is an even more fragile property; although many full Euler polars are convex, there is no simple mathematical condition to ensure it, other than the restrictive assumption of constant heat capacity (polytropic). On the other hand \cite{elling-idealpolar} has shown using only standard assumptions that critical points are unique and that there are unique weak and strong shocks, so again these properties are more robust than polar convexity. 

For the simplest non-ideal case, the van der Waals eos, \cite{elling-waalspolar} found numerical counterexamples in physically reasonable regions of the thermodynamic phase plane. Only in the barotropic Euler case convex eos \emph{happens} to imply strictly convex shock polar (\cite{elling-isenpolar}; see also \cite{preiswerk-phdthesis1938}). Earlier discussion of shock polars can be found in \cite{teshukov-polar,henderson-menikoff,fowles-jfm1981}; see \cite{bethe,weyl-shock-waves} for essential background on normal shocks for general eos.

\section{Full Euler equations of state}

We consider some standard assumptions about the eos $\epm=\epm(\idens,\spm)$. We only consider smooth functions $\epm$ that satisfy all of the following:\\
\begin{enumerate}
\item
  The domain of $\epm$ contains only $(\idens,\spm)$ with $\idens>0$. 
\item
  Positive temperature $\temp=\epm_\spm$ and pressure $\pp=-\epm_\idens$.
\item
  The matrix of second derivatives
  \[ \begin{bmatrix}
    \epm_{\spm\spm} & \epm_{\idens\spm} \\ 
    \epm_{\idens\spm} & \epm_{\idens\idens} 
  \end{bmatrix} \] 
  is positive definite, which means that each state $(\idens,\spm)$ is thermodynamically stable. In particular $\epm_{\idens\idens}>0$ implies $\csnd$ is well-defined.
  This condition is violated in certain regions of van der Waals eos unless modifications are made to produce evaporation-condensation equilibrium.
\item
  The eos is \defm{convex}, with ``convex'' referring to $\idens\mapsto\pp(\idens,\spm)$:
  \[ \pp_{\idens\idens} > 0 . \]
  Due to $\pp_\idens=-(\dens\csnd)^2$ the condition is equivalent to
  \[ \dens(\frac{\partial(\csnd^2) }{\partial\dens})_\spm > - 2\csnd^2 . \myeqlab{eq:convexeos-cc}\]
  When this condition is violated, even normal steady shocks of weak strength are not determined uniquely by the upstream state,
  allowing polars to enclose a nonconvex or even disconnected area; 
  besides, the physically reasonable solution to corner reflections and many other problems may not be the compressive shock but rather a fan or some shock-fan combination,
  so that the meaning of the unmodified shock polar is dubious without modification.
\end{enumerate}

If the eos additionally satisfies the ideal gas law
\[ \frac{\pp\idens}{\temp} = \Rspec = \text{const} , \myeqlab{eq:idealgaslaw}\]
then it satisfies Joule's second law: $\epm=\epm(\temp)$ is a function of $\temp$ alone. 
This can be seen by substituting $\temp=\epm_\spm$ and $\pp=-\epm_\idens$ into \myeqref{eq:idealgaslaw}, which yields a first-order PDE
\[ 0 = \Rspec\epm_\spm + \idens\epm_\idens \]
with general solution
\[ \epm = \epm(\subeq{\frac{\spm-\spmu}{\Rspec}-\ln\frac{\idens}{\idensu}}{\xsv}) . \]
Note
\[ \Rspec\temp = \Rspec\epm_\spm = \epm_\xsv \]
and 
\[ \epm_{\xsv\xsv} = \Rspec^2\epm_{\spm\spm} > 0 \]
by assumption 3 above, so that $\xsv\mapsto \epm_{\xsv}=\Rspec\temp$ is a strictly increasing function; we can invert it and obtain $\xsv$ as a function of $\temp$,
so that $\epm=\epm(\temp)$.

Restating the definition of $\temp,\pp$ as the fundamental relation
\[ d\epm = \temp d\spm - \pp d\idens \topref{eq:idealgaslaw}{=} \temp d\spm - \Rspec \temp \frac{d\idens}{\idens} , \]
we derive
\[ \frac{\spm-\spmu}{\Rspec} = \ln\frac{\idens}{\idensu} + \int \frac1{\Rspec\temp} d\epm . \myeqlab{eq:svT-ideal} \] 

Enthalpy $\hpm=\epm+\pp\idens=\epm+\Rspec\temp$ also becomes a function of $\temp$ alone, as does 
sound speed $\csnd^2$:
\[ (\dens\csnd)^2 &= \epm_{\idens\idens} = (-\epm_\xsv/\idens)_{\idens} = (\epm_{\xsv\xsv}+\epm_\xsv)/\idens^2
\\&= (\Rspec^2\epm_{\spm\spm}+\Rspec\epm_\spm)\dens^2 = \Rspec\epm_\spm(\frac{\Rspec\temp_\spm}{\epm_\spm}+1)\dens^2 = \Rspec\temp(\Rspec\temp_\epm+1)\dens^2 
\\\impl\quad \csnd^2 &= \Rspec\temp(\Rspec\temp_\epm+1) . \myeqlab{eq:cct} \]

Monotone sound speed corresponds to
\[ 0 &\leq (\csnd^2)_\temp = (\csnd^2)_\epm \subeq{\epm_\temp}{>0}
\\\eqv\quad 0 &\leq (\Rspec\temp(\Rspec\temp_\epm+1))_\epm = \Rspec\temp \Rspec\temp_{\epm\epm} + \Rspec\temp_\epm(\Rspec\temp_\epm+1) . \]
Using $\temp_{\epm\epm}=-\epm_{\temp\temp}/\epm_\temp^3$ we transform to
\[ \epm_{\temp\temp} &\leq \frac{\epm_\temp(\epm_\temp+\Rspec)}{\Rspec\temp}. \]
The borderline case of $=$ has general solution 
\[ \epm_\temp = \frac{\Rspec}{C/\Rspec\temp-1}  \myeqlab{eq:borderline-mono} \]
for constants $C$.

Convex eos is a looser condition: 
\[ -2\csnd^2 \topref{eq:convexeos-cc}{<} -\idens(\frac{\partial(\csnd^2)}{\partial\idens})_\spm = -\idens\frac{\partial(\csnd^2)}{\partial\epm} \subeq{(\frac{\partial\epm}{\partial\idens})_\spm}{=-\pp=-\Rspec\temp/\idens} = (\csnd^2)_\epm \Rspec \temp . \]
After filling in \myeqref{eq:cct} some calculation produces
\[ \epm_{\temp\temp} < \frac{ \epm_\temp(\epm_\temp+\Rspec)(2\epm_\temp+\Rspec) }{ \Rspec^2\temp } . \]
The borderline case with $<$ replaced by $=$ has general solution
\[ \epm_\temp = \Rspec\frac{(1-4C\Rspec\temp)^{-\frac12}-1}{2} .  \myeqlab{eq:borderline-convex} \]

Polytropic eos have constant heat capacities $\cvspec=\epm_\temp$ and therefore $\cpspec=\hpm_\temp=\Rspec+\cvspec$.
Heat capacity is nearly constant over wide temperature ranges for noble gas, narrow ranges for oxygen, nitrogen, hydrogen and other symmetric diatomic molecules;
for carbon dioxide, hydrocarbons and more complex molecules polytropic constant $\cvspec$ can be regarded as a useful first approximation.
\[ \gisen=\cpspec/\cvspec \]
is called ratio of heats, or adiabatic/isentropic exponent. We will later use
\[ \hpm = \cpspec\temp = \subeq{\frac{\cpspec}{\Rspec}}{=\gisen/(\gisen-1)}\pp\idens . \]
We also simplify
\[ \csnd^2 \topref{eq:cct}{=} \Rspec\temp(\frac{\Rspec}{\cvspec}+1) = \gisen\Rspec\temp = \gisen\pp\idens . \myeqlab{eq:ccpoly} \] 

\myeqref{eq:svT-ideal} reduces to
\[ \frac{\spm - \spmu}{\Rspec} = \ln\frac{\idens}{\idensu} + \frac{\cvspec}{\Rspec} \ln\frac{\epm}{\epmu}
= \ln\frac{\idens}{\idensu} + \frac1{\gisen-1} \subeq{\ln\frac{\epm}{\epmu}}{=\ln(\temp/\tempu)} . \]
For models like potential flow or isentropic Euler that use constant entropy $\spm$, we substitute $\temp=\pp\idens/\Rspec$ again to obtain 
\[ \pp = \ppu (\frac{\idens}{\idensu})^{-\gisen} . \]

\section{Initial steps common to all models}

Many properties of shock polars can be shown by some brute-force approach, but those lead through rather complicated arguments and calculations; worse, they produce theorems with assumptions stronger than necessary. The key to minimal assumptions and short proofs is a careful choice of methods. We start with a familiar observation about density that all models have in common.

\begin{figure}
  \parbox{.45\linewidth}{%
    \centerline{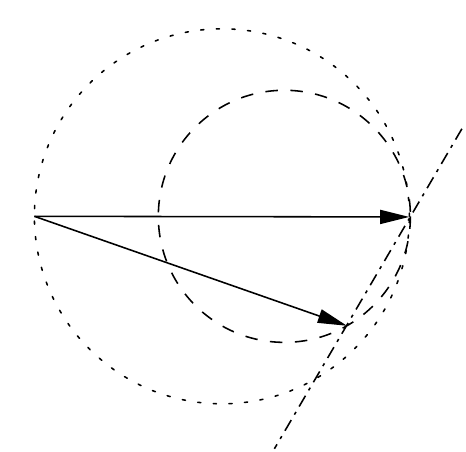}
    \caption{Level sets of density $\dens=\dens(\vv)$ are circles with diameter from $\vvu\idens/\idensu$ to $\vvu$.}
    \label{fig:rhocircle}
  }\hfil\parbox{.45\linewidth}{%
    \centerline{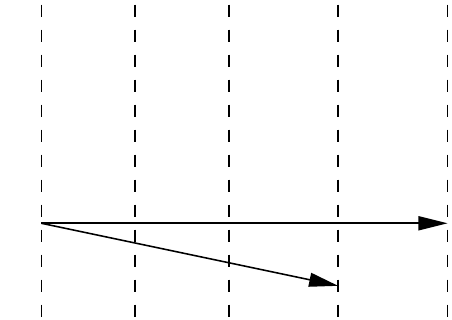}
    \caption{Euler case: level sets of pressure $\pp=\pp(\vv)$ are vertical lines, with pressure increasing from right to left.}
    \label{fig:plines}
  }
\end{figure}

Both full Euler and full potential flow are rotationally symmetric models, so we may rotate coordinates to make $\vvu$ horizontal rightward, i.e.\ $\vxu>0$ and $\vyu=0$.
Then $\vvu/|\vvu|=(1,0)$; we abbreviate
\[ (\xiv,\etav)= \vv/|\vvu| . \]

By \myeqref{eq:full-vt} velocity jump $\vvu-\vv$ is normal to the shock, and since $\vvu-\vv$ points downstream for admissible shocks,
\[ \nn = \frac{\vvu-\vv}{|\vvu-\vv|} \]
is the unique downstream normal. We use this to simplify mass conservation \myeqref{eq:full-jn}:
\[ 1 = \frac{\densu\vvu\dotp\nn}{\dens\vv\dotp\nn} = \frac{\densu\vvu\dotp(\vvu-\vv)}{\dens\vv\dotp(\vvu-\vv)} . \]
After some calculation this reduces to the convenient formula
\[ \frac{\dens}{\densu} = \frac{1-\xiv}{\xiv(1-\xiv)-\etav^2} .\myeqlab{eq:dens-xy} \] 
Alternatively we can solve the equation for
\[ \etav^2 = (\xiv-\frac{\idens}{\idensu})(1-\xiv) . \myeqlab{eq:yyx}\]
For constant $\idens$ this equation describes $(\xiv,\etav)$ on a circle with diameter the line segment from $(\idens/\idensu,0)$ to $(1,0)$ (fig.\ \myref{fig:rhocircle}).
The interior of the circle for $\idens=0$ is called the \defm{compressive disk}; its interior contains exactly those $(\xiv,\etav)$ that correspond to compressive shocks.
Adjacent are points with nonsensical negative or undefined $\idens$, while the $\xiv>1$ halfplane has $\idens$ that are finite and positive, but larger than $\idens_0$, i.e.\ non-compressive shocks.

\section{Full Euler}

The normal momentum conservation \myeqref{eq:full-momn}
\[ \pp + \dens(\vv\dotp\nn)^2 = \ppu + \densu(\vvu\dotp\nn)^2 \]
reduces, using the mass relation \myeqref{eq:full-jn} $\dens\vv\dotp\nn=\densu\vvu\dotp\nn$, to
\[ \pp = \ppu + \densu\vvu\dotp\nn (\vvu-\vv)\dotp\nn . \]
Since velocity jump $\vvu-\vv$ is normal to the shock, we may again use $\nn=(\vvu-\vv)/|\vvu-\vv|$, to simplify:
\[ \pp &= \ppu + \densu\vvu\dotp\frac{\vvu-\vv}{|\vvu-\vv|} (\vvu-\vv)\dotp\frac{\vvu-\vv}{|\vvu-\vv|} = \ppu + \densu\vvu(\vvu-\vv) . \]
Written with $(\xiv,\etav)$ this is
\[ \frac{\pp-\ppu}{\densu|\vvu|^2} &= 1-\xiv \myeqlab{eq:pj1x}. \]
So the pressure is essentially a horizontal coordinate in the $\vv$ polar diagram (see fig.\ \myref{fig:plines}), a fact already noted and used by Busemann and other classical authors.

\subsection{Polytropic full Euler}

We use the Bernoulli relation \myeqref{eq:full-Hpm}, \myeqref{eq:Hpm}:
\[ 0 = 2\jmp\hpm + \jmp{|\vv|^2} \]
After division by $|\vvu|^2$, some calculation shows
\[ 0 = 2\frac{\hpm-\hpmu}{|\vvu|^2} + \xiv^2+\etav^2 - 1 . \]
For polytropic eos we can combine this with 
\[ 2\hpm = 2\cpspec\temp = \frac{2\cpspec}{\Rspec}\pp\idens = \frac{2\gisen}{\gisen-1} \pp\idens , \]
as well as \myeqref{eq:pj1x} $\jmp\pp/\densu|\vvu|^2=1-\xiv$ and \myeqref{eq:dens-xy} $\idens/\idensu=(\xiv(1-\xiv)-\etav^2)/(1-\xiv)$, to get a relation linear in $\etav^2$ with coefficients polynomial in $\xiv$:
\[ 0
= \frac{2\gisen}{\gisen-1} \Big( (\frac{\ppu}{\densu|\vvu|^2}+1-\xiv)\frac{\xiv(1-\xiv)-\etav^2}{1-\xiv} - \frac{\ppu}{\densu|\vvu|^2} \Big) + \xiv^2+\eta^2-1 \] 
After some calculation the solution is obtained as
\[ \etav^2 = (1-\xiv)^2 \frac{\xiv-\xinor}{\ximax-\xiv} \myeqlab{eq:poly-euler-polar} \]
where 
\[ \xinor = \frac{2\gisen\frac{\ppu}{\densu|\vvu|^2}+\gisen-1}{\gisen+1}   \csep   \ximax = 1+\frac{2\gisen\frac{\ppu}{\densu|\vvu|^2}}{\gisen+1} . \]
Using
\[ \frac{\gisen\ppu}{\densu|\vvu|^2} \topref{eq:ccpoly}{=} \frac{\gisen}{\densu|\vvu|^2}\frac{\csndu^2}{\gisen\idensu} = \frac1{\Machu^2} , \]
this simplifies to 
\[ \xinor = \frac{2\Machu^{-2}+\gisen-1}{\gisen+1}   \csep   \ximax = 1+\frac{2\Machu^{-2}}{\gisen+1} . \]
\myeqref{eq:poly-euler-polar} is classical; in textbooks it is usually derived via the Prandtl relation in a longer calculation (see \cite[section 121]{courant-friedrichs}
or \cite[p.\ 11]{ferrari-tricomi}).

Consider the formula for the upper half-polar:
\[ \etav \topref{eq:poly-euler-polar}{=} |\xiv-1||\xiv-\xinor|^{1/2}|\xiv-\ximax|^{-1/2} . \myeqlab{eq:eta-upperhalf}\] 
We note that for $\gisen\in\boi{-1}{\infty}$ and $\Machu\in\boi1\infty$ we have
\[ \xinor &= \frac{\gisen-1+2\Machu^{-2}}{\gisen+1} \in \boi{-\infty}{1},
\\ \ximax &= 1 + \frac{2\Machu^{-2}}{\gisen+1} \in \boi{1}{\infty}. \]

Generally $\etav=\prod_j|\xiv-a_j|^{\alpha_j}$ has derivative
\[ \etav_\xiv = \etav l \] 
where $l=\etav_\xiv/η=(\ln \etav)_\xiv$ is the logarithmic derivative
\[ l = \sum_j\frac{\alpha_j}{\xiv-a_j} . \]
\[ \etav_{\xiv\xiv} = \etav_\xiv l+\etav l_\xiv = \etav(l²+l_\xiv) = \etav \Big( (\sum\frac{\alpha_j}{\xiv-a_j})^2 - \sum \frac{\alpha_j}{(\xiv-a_j)^2} \Big).  \myeqlab{eq:logderder}\]
$η>0$ by definition when $\xiv\in\boi{\xinor}{1}\subset\boi{\xinor}{\ximax}$; in our case the big parenthesis of \myeqref{eq:logderder} is (cf.\ \myeqref{eq:eta-upperhalf})
\[ &=
\Big( \subeq{\frac{1}{\xiv-1}}{\eqdef-U} + \half\subeq{\frac1{\xiv-\xinor}}{\eqdef n} + (-\half)\subeq{\frac1{\xiv-\ximax}}{\eqdef -m} \Big)^2 - \subeq{\frac1{(\xiv-1)^2}}{=U^2} - \half\subeq{\frac1{(\xiv-\xinor)^2}}{=n^2} - (-\half)\subeq{\frac1{(\xiv-\ximax)^2}}{=m^2}
\intertext{(where we chose signs so that $U,m,n>0$ for $\xinor<\xiv<1<\ximax$)}
&=
U^2 - U(n+m) + \frac14 (n+m)^2 \qquad - U^2 - \frac12 (n^2-m^2)
\\&=
(n+m)\frac14\big(-4U+(n+m)-2(n-m)\big)
\\&=
(n+m)\frac14\big(-4U-n+3m\big).
\myeqlab{eq:etaxx-eta}
\]
$n+m$ is positive. For the last parenthesis we use $\xiv<1<\ximax$ to argue $0<1-\xiv<\ximax-\xiv$, so $U=1/(1-\xiv)>1/(\ximax-\xiv)=m$, so rearranging said last parenthesis to
\[ -3(\subeq{U-m}{>0}) - U - n \]
we confirm it is negative. 
Altogether we have shown $\etav_{\xiv\xiv}<0$, so $η$ has concave graph over $\xiv∈\boi{\xinor}{1}$. 

We have not dealt with the $\xiv=\xinor$ (normal) or $\xiv=1$ (vanishing) points. The vanishing point is a strictly convex corner where the upper half-polar is well-known to have slope $-\tan\arcsin(1/\Machu)\in\boi{-\infty}{0}$, with $\arcsin(1/\Machu)$ the \defm{Mach angle} that acoustic fronts form to the upstream velocity. At the normal point we can use $\xiv>\xinor$ for any other point to argue strict convexity. 

We recover a result that is essentially classically known, but often not stated explicitly or for the same parameter range:
\begin{theorem}
  For any $-1<\gisen<\infty$ and any $1<\Machu<\infty$ the polytropic full Euler compressive velocity-plane shock polar is a strictly convex curve.
\end{theorem}
(Of course ``strictly convex'' has the meaning consistent with the literature, that as we move from right to left (larger to smaller $\vx$) the slope of the upper halfpolar is strictly increasing, with symmetric result for the lower halfpolar.)

We caution that for $\gisen<1$ and large $\Machu$ the polar defined by \myeqref{eq:poly-euler-polar} can have unusual features, namely leaving the compressive circle. The formula \myeqref{eq:poly-euler-polar} happens to have, for polytropic eos, a mathematical extension that has no physical meaning; for most other eos some singular behaviour at the boundary of the compressive circle would prevent extension.

\subsection{Ideal non-polytropic full Euler}

\begin{figure}
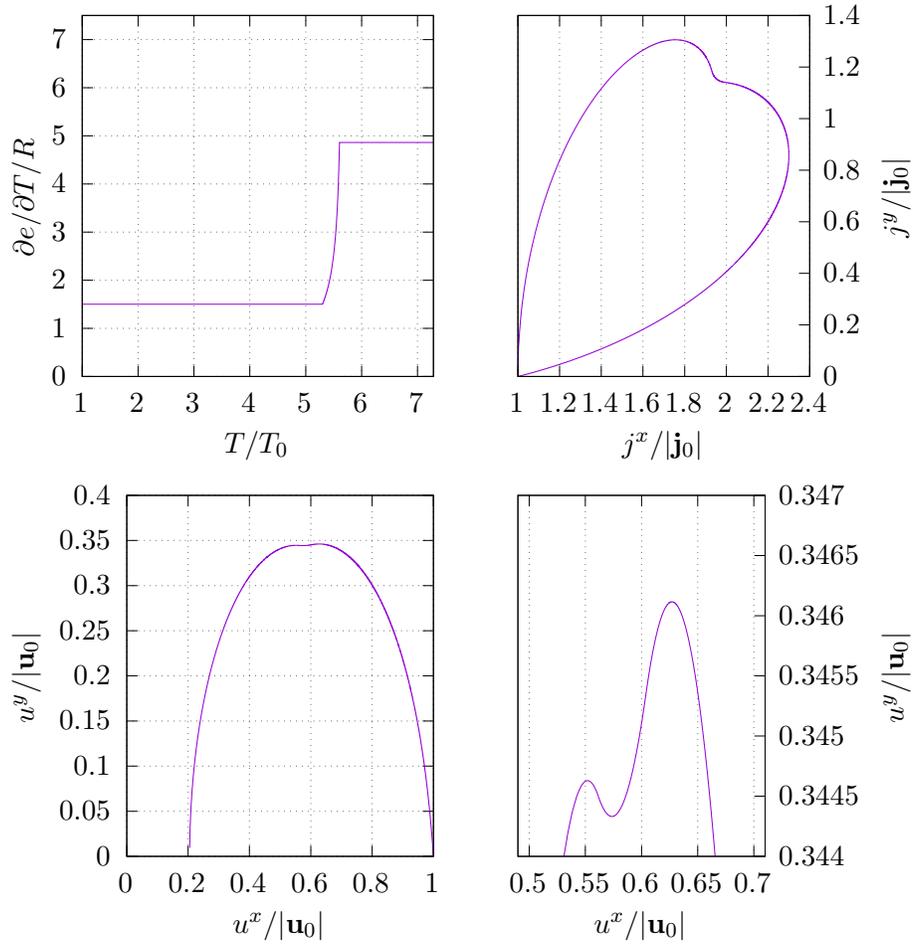

\input{idealeuler-nonconvex-vv-polar-cv.xxx}\input{idealeuler-nonconvex-jj-polar.xxx} \\
\input{idealeuler-nonconvex-vv-polar.xxx}\input{idealeuler-nonconvex-vv-polar-detail.xxx}
\caption{Ideal Euler with $\partial\epm/\partial\temp$ as shown in the upper left diagram and $\Machu=5.0$ produces a clearly non-convex $\jj$ polar (upper right) and a slightly non-convex $\vv$ polar (lower left, detail in lower right). The equation of state is borderline convex in the non-constant region of $\partial\epm/\partial\temp$ for $\temp/\tempu$ (from $5.3$ to about $5.6$).}
\label{fig:idealpolarnonconvex}
\end{figure}

For full Euler with an ideal but non-polytropic eos, we revisit momentum conservation \myeqref{eq:full-momn}:
\[ \jnq = \frac{\jmp\pp}{-\jmp{\idens}} . \myeqlab{eq:pvjnq} \]
We combine this with the Bernoulli relation \myeqref{eq:full-Hpm},\myeqref{eq:Hpm}: 
\[ 0 &= 2\jmp{\hpm} +\jmp{|\vv|²} \overset{\jmp{\vt}=0}{=} 2\jmp{\hpm} + \jmp{\vnq}
\\&= 2\jmp{\hpm} + \jmp{\idens^2\jnq} \overset{\jmp{\jn}=0}{=} 2\jmp{\hpm} + \jnq\jmp{\idens^2} . \]
Using $\jmp{ab}=\avg{a}\jmp b+\jmp{a}\avg b$, where $\avg{a}=(a+a_0)/2$ is the average of up- and downstream quantities, so that $\jmp{a^2}=2\jmp{a}\avg{a}$, the last relation turns into 
\[ 0 &= \jmp{\hpm} + \jnq \avg\idens \jmp\idens
\topref{eq:pvjnq}{=} \jmp{\hpm} - \avg\idens \jmp{\pp}. \]
This is the well-known \defm{Hugoniot relation}, distinguished because it relates purely thermodynamic quantities without velocity involved.
The same name is applied to a common variant, derived from the substitution $\hpm=\epm+\idens\pp$:
\[ 0 &= \jmp{\epm} + \subeq{\jmp{\idens\pp}}{=\jmp{\idens}\avg\pp+\jmp\pp\avg\idens} - \avg\idens\jmp\pp = \jmp\epm + \avg\pp \jmp\idens . \]

For ideal gas $\epm=\epm(\temp)$. Besides, the ideal gas law \myeqref{eq:idealgas} gives $\idens=\Rspec\temp/\pp$, so that 
\[ 0 = \jmp{\epm(\temp)} + \avg\pp \jmp{\Rspec\temp/\pp} \]
which is essentially a quadratic equation for $\pp$, when $\temp,\tempu,\ppu$ are given. After some calculation the solution turns out to be 
\[ \frac{\pp}{\ppu} = \frac{\jmp{\epm/\Rspec+\half\temp}}{\tempu} + \sqrt{\Big(\frac{\jmp{\epm/\Rspec+\half\temp}}{\tempu}\Big)^2+\frac{\temp}{\tempu}} ; \]
the second solution has $-\sqrt{...}$, but for positive temperatures $\temp,\tempu$ for realistic eos, this solution always has negative pressure $\pp$, so we ignore it.

Hence, given $\temp$, there is an explicit formula for $\pp$. By \myeqref{eq:pj1x} $\pp$ gives $\xiv$, while the ideal gas law provides $\idens=\Rspec\temp/\pp$. Finally, conservation of mass in the form \myeqref{eq:yyx} $\etav^2=(\xiv-\idens)(1-\xiv)$ gives $\etav^2$. Therefore we have a chain of explicit formulas for calculating the shock polar $\temp\mapsto(\xiv,\etav)$ parametrized by $\temp$.

Trying to parametrize by $\xiv$, as in the polytropic case, is generally difficult, since $\epm=\epm(\temp)$ can have a rather arbitrary function $\epm$ that is generally hard to invert explicitly, except in special cases like the polytropic one where $\epm=\cvspec\temp$ for constant $\cvspec$.

Given the explicit formulas it is straightforward to plot graphs of the shock polar for any given function $\epm(\temp)$.
While \cite{elling-idealpolar} proved, by rather long calculations and arguments, that critical points are always unique and subsonic for ideal eos satisfying standard assumptions (in particular that the eos is convex),
the resulting shock polars need \emph{not} be convex under the same assumptions.
An example is shown in figure \myref{fig:idealpolarnonconvex}. There, $\epm(\temp)$ is linear over two ranges of temperatures $\temp$, but in between it is chosen to give a borderline convex eos; the non-horizontal segment in figure \myref{fig:idealpolarnonconvex} top left is an exact solution of this condition expressed as an ODE constraining $\epm(\temp)$, as in \myeqref{eq:borderline-convex}. 

The $\vv$ polar is slightly non-convex (fig.\ \myref{fig:idealpolarnonconvex} bottom left and detail bottom right); the $\jj$ polar for mass flux $\jj$ is more clearly non-convex (top right).

It is not difficult to give rigorous proofs of these or later counterexamples, most trivially by taking two derivatives to calculate curvature, then evaluating the formula at a selected point using exact or interval arithmetic, or more elegantly by taking certain asymptotic limits, but we leave it to interested readers since the numerical plots are quite clear.

Although the rise of heat capacity shown in fig.\ \myref{fig:idealpolarnonconvex} top left is quite abrupt, steeper than in most realistic fluids, it does demonstrate that there is no easy theorem that gives convexity of $\vv$ polars, unlike the main theorem in \cite{elling-idealpolar} which proved uniqueness of critical shocks assuming only convex eos. 

Non-convexity is still possible when restricting the eos to those with monotone sound speed, meaning $\csnd$ is a nondecreasing function of $\temp$.
For example for $\Machu=45$, an eos $\epm(\temp)$ that is polytropic with $\fod=3$ for $\temp<\temp_1=10.3\tempu$ is followed by an ideal constant-$\csnd$ branch (as in \myeqref{eq:borderline-mono}) for $\temp_1<\temp<\temp_2$ with $\temp_2=1.05\temp_1$, then followed by another polytropic branch on the $\temp>\temp_2$ side, with free constants chosen to make $\epm_\temp$ continuous across the $\temp=\temp_1,\temp_2$ knots. 
The resulting $\vv$ polar is non-convex at $\temp$ slightly larger than $\temp_1$. 
The non-convexity is more slight in these cases and is not represented graphically here.

Primarily these examples show that there appears to be no reasonably simple and correct theorem giving convex polars for a wide range of non-polytropic Euler polars, even after restriction to the ideal case. Consequently the assumption of convexity is unsatisfying and should be avoided after moving beyond polytropic eos.

\section{Full potential flow}

\begin{figure}
  \centerline{\input{polarpotf-vv-nonconvex.xxx}}
  \input{polarpotf-vv-nonconvex-detail.xxx}
  \input{polarpotf-vv-nonconvex-c.xxx}
  \caption{Non-convex $\vv$ polar for full compressible potential flow with $\gisen=-0.75$ for $\idens$ above $0.019\idensu$, $\gisen=5/3$ below}
  \label{fig:potf-vv-nonconvex}
\end{figure}

In the compressible full potential flow case, $\spm$ is treated as constant, so that the eos is specified by a scalar function of a single variable, for example $\pp=\pp(\dens)$. This then produces $\hpm=\hpm(\dens)$ by the relation
\[ \hpm_\dens = \frac{\pp_\dens}{\dens} = \frac{\csnd^2}{\dens} ; \]
we assume of course that $\csnd^2>0$ so that sound speed $\csnd$ is well-defined (we say the eos is \defm{hyperbolic}). Fixing a constant $\spm$ for polytropic eos produces
\[ \pp=\ppu(\dens/\densu)^\gisen . \]
The potential flow shock relations \myeqref{eq:potf-jn} to \myeqref{eq:potf-Hpm} do not involve $\pp$ at all, so that it is mathematically sufficient to specify some $\hpm(\dens)$ with positive derivative. 

Shock relations $\jmp{\jj\dotp\nn}=0$ and $\jmp{\vv\dotp\ts}=0$ are still available for potential flow and still imply \myeqref{eq:yyx}:
\[ \etav^2 = (\xiv-\frac{\idens}{\idensu})(1-\xiv) . \myeqlab{eq:potf-yyx}\] 
However, since the normal momentum relation \myeqref{eq:full-momn} is no longer available for potential flow, neither is the $\pp$-$\xiv$ relation in \myeqref{eq:pj1x}.
$\idens$ is instead defined implicitly by the Bernoulli relation \myeqref{eq:potf-Hpm} (note \myeqref{eq:Hpm}):
\[ - 2\frac{\jmp\hpm}{|\vvu|^2} = \xiv^2+\etav^2 - 1 . \]
Again $\etav^2$ appears linearly, so we can eliminate it using \myeqref{eq:potf-yyx}:
\[ -2\frac{\jmp\hpm}{|\vvu|^2} = \xiv^2 + (\xiv-\frac{\idens}{\idensu})(1-\xiv) - 1 = (\xiv-1)(1+\frac{\idens}{\idensu}). \]
This is now easily solved for $\xiv$ as a function of $\idens$:
\[ \xiv = 1 + \frac{-2(\hpm-\hpmu)/|\vvu|^2}{1+\idens/\idensu} . \myeqlab{eq:x-idens} \]
\myeqref{eq:potf-yyx} and \myeqref{eq:x-idens} combined give a perfectly explicit formula for the shock polar parametrized by $\idens$ for any function $\hpm(\idens)$.

However, to give an explicit formula parametrized by $\xiv$, for example by inverting \myeqref{eq:x-idens} to $\idens$ a function of $\xiv$, is generally not possible if $\hpm=\hpm(\idens)$ is permitted to be a rather arbitrary function. This is similar to the full Euler situation when passing from polytropic to general eos. 

But for our purposes we only need existence of an inverse, not an explicit formula. 
Take the $\idens$ derivative of the right-hand side of \myeqref{eq:x-idens} at some $\idens<\idensu=1$:
\[  0 + \frac1{|\vvu|^2} \Big( \frac{-2\hpm_\idens}{1+\idens/\idensu} + \frac{2(\hpm-\hpmu)/\idensu}{(1+\idens/\idensu)^2} \Big) . \]
This is positive because $\idens\hpm_\idens=-\dens\hpm_\dens=-\csnd^2<0$ so that $\hpm_\idens<0$ and thus $\hpm>\hpmu$ since $\idens<\idensu$.
Therefore the inverse function theorem can be used to invert \myeqref{eq:x-idens}, showing that $\xiv\mapsto\idens$ is a well-defined smooth strictly increasing function on any interval of $\idens<\idensu$ in the domain of definition of $\hpm$.

Consider some $\idens<\idensu=1$ so that \myeqref{eq:x-idens} yields $\xiv\in\boi{\idens/\idensu}{1}$, so that $\etav>0$ is well-defined by \myeqref{eq:potf-yyx}, with $(\xiv,\etav)$ inside the compressive circle.
To check strict concavity of $\xiv\mapsto\etav$ we take two derivatives of
\[ \etav = \sqrt{ab} \quad\text{where}\quad a = \xiv-\frac{\idens}{\idensu} \csep b = 1-\xiv. \]
Result:
\[ \etav_\xiv &= \frac{a_\xiv b+ab_\xiv}{2\etav}, 
\\ \etav_{\xiv\xiv}
&= \frac{a_{\xiv\xiv}b+2a_\xiv b_{\xiv}+ab_{\xiv\xiv}}{2\etav} - \frac{(a_{\xiv}b+ab_{\xiv})^2}{4\etav^3} 
\\&= \frac{2a_{\xiv\xiv}ab^2+4a_{\xiv}b_{\xiv}ab+2a^2bb_{\xiv\xiv}}{4\etav^3} - \frac{a_{\xiv}^2b^2+2aba_{\xiv}b_{\xiv}+a^2b_{\xiv}^2}{4\etav^3} 
\\&= \frac{2a_{\xiv\xiv}ab^2+2a^2bb_{\xiv\xiv}}{4\etav^3} - \frac{(ba_{\xiv}-ab_{\xiv})^2}{4\etav^3} . \]
Using $a_{\xiv\xiv}=-\idens_{\xiv\xiv}/\idensu$ and $b_{\xiv\xiv}=0$ we get
\[ \etav_{\xiv\xiv}
&=  \frac{-2\frac{\idens_{\xiv\xiv}}{\idensu}(\xiv-\frac{\idens}{\idensu})(1-\xiv)^2}{4\etav^3} - \frac{(ba_\xiv-ab_\xiv)^2}{4\etav^3},  \]
so if $\idens_{\xiv\xiv}>0$, then due to $\idens/\idensu<\xiv<1$ we find $\etav$ is a strictly concave function.

To estimate $\idens_{\xiv\xiv}$, we take $d/d\xiv$ derivatives of \myeqref{eq:x-idens} in the form
\[ 0 = 2\frac{\hpm-\hpmu}{|\vvu|^2} + (\xiv-1)(1+\frac{\idens}{\idensu}) ,\]
giving 
\[ 0 = 2\frac{\hpm_\idens\idens_\xiv}{|\vvu|^2} + 1 + \frac{\idens}{\idensu} + \frac{\idens_\xiv}{\idensu}(\xiv-1)
\\\impl\quad \frac{\idens_\xiv}{\idensu}
= \frac{1+\idens/\idensu}{1-\xiv-2\hpm_\idens\idensu/|\vvu|^2} \]
which is positive because $\hpm_{\idens}=-\csnd^2/\idens<0$ and $\xiv<1$ so that the denominator is positive. Using this we also find that
\[ \frac{\idens_{\xiv\xiv}}{\idensu}
&= \subeq{\frac{\idens_\xiv/\idensu}{1-\xiv-2\hpm_\idens\idensu/|\vvu|^2}}{>0}+\frac{(1+\idens/\idensu)(1+2\hpm_{\idens\idens}\idensu\idens_\xiv/|\vvu|^2)}{(1-\xiv-2\hpm_\idens\idensu/|\vvu|^2)^2} \]
is positive if $\hpm_{\idens\idens}\geq 0$.

We rephrase this condition in terms of sound speed:
\[ \hpm_{\idens\idens}
=
\partial_\idens(-\frac{\csnd^2}{\idens})
=
\frac{-(\csnd^2)_\idens}{\idens}+\frac{\csnd^2}{\idens^2}
=
\frac{(\csnd^2)_\dens}{\idens^3}+\frac{\csnd^2}{\idens^2}
=
\dens^2\Big( \dens (\csnd^2)_\dens + \csnd^2 \Big).
\]
Clearly $\hpm_{\idens\idens}\geq 0$ \emph{is implied} by \defm{monotone sound speed}, meaning
\[ \csnd_\dens > 0 . \]
On the other hand $\hpm_{\idens\idens}\geq 0$ \emph{implies}
\[ \dens(\csnd^2)_\dens \geq - \csnd^2 > -2 \csnd^2 , \]
which is the condition of \defm{convex eos} in the form \myeqref{eq:convexeos-cc}. In short, 
\[ \text{monotone sound speed} \impl \text{convex $\hpm$} \impl \text{convex eos}. \]

The normal and vanishing point are treated by arguments analogous to the full Euler case. 

In summary:
\begin{theorem}
  For a hyperbolic eos
  the compressive part of the potential flow velocity plane shock polar is strictly convex
  if enthalpy $\hpm$ is a convex function of $\idens$.
  This holds in particular if sound speed $\csnd$ is a nondecreasing function of $\dens$. 
\end{theorem}
(Again, ``strictly convex'' has the meaning consistent with the literature, that as we move from right to left the slope of the upper halfpolar is strictly increasing.)

We examine the $\gisen$-law case: omitting irrelevant additive constants we have
\[ \pp = \ppu(\frac{\dens}{\densu})^{\gisen} , \]
giving
\[ \csnd^2 &= \pp_\dens = \frac{\gisen\ppu}{\densu} (\frac{\dens}{\densu})^{\gisen-1}  , \]
showing we need to restrict $\ppu$ to $\sign\ppu=\sign\gisen$ for hyperbolicity (by adding a constant, pressure can be made positive on some interval, if desired for physical reasons, but mathematically there is no difference).
$\hpm_\dens=\csnd^2/\dens$ fixes, up to an additive constant, 
\[ \hpm &= \frac{\gisen\ppu}{(\gisen-1)\densu} \subeq{ (\frac{\dens}{\densu})^{\gisen-1} }{ = (\idens/\idensu)^{1-\gisen} }. \]
Then
\[ \hpm_\idens = - \gisen\ppu (\frac{\idens}{\idensu})^{-\gisen} < 0 , 
\\ \hpm_{\idens\idens} = \gisen^2\ppu (\frac{\idens}{\idensu})^{-\gisen-1} , \]
which is positive for $\gisen>0$ so that $\ppu>0$ by assumption, but not for $\gisen<0$ where $\ppu<0$. 
So the theorem proves strict convexity of the polar for $\gisen>0$.

Counterexamples can be found for cases of $-1<\gisen<0$. Figure \myref{fig:potf-vv-nonconvex} displays one of them. In the figure the eos was switched from $\gisen=-0.75$ for $\idens$ larger than $0.019\idensu$ to $\gisen=5/3$ at smaller values, to make the polar complete. 

For compressible potential flow, strictly convex shock polars are therefore a bit more common than for full Euler with ideal eos. The theorem is a significant extension of prior work, namely \cite[Theorem 1]{elling-sonic-potf}, which only gave an ad-hoc proof in the special case of $\gisen>1$ polytropic eos. The old proof does not generalize well because it assumes monotonicity of downstream normal velocity and other quantities along the polar. These are convenient but fragile properties that do not generalize well beyond the polytropic special cases.

\section{Conclusions}

The convexity theorems do not and cannot cover the whole range of convex eos, so both in the full Euler and potential flow case it is preferable to avoid the property, except in a first approach. 
After all, it was shown in \cite{elling-polarpotf} that for any convex eos critical points are unique and subsonic, which is for many purposes the property actually needed. 
Besides, uniqueness of critical, weak and strong shocks is about how often certain curves intersect, which is a coordinate-independent property,
unlike convexity, which could change from a $\vv$ polar diagram to the familiar $\pp$-$\turn$ plane, to the $\jj$ plane etc.

\section*{Acknowledgement}

This research was partially supported by Taiwan MOST Grant No.\ 110-2115-M-001-005-MY3.

\bibliographystyle{amsalpha}

\begin{thebibliography}{Wey49}

\bibitem[Bet42]{bethe}
{H.} Bethe, \emph{On the theory of shock waves for an arbitrary equation of
  state}, Tech. Report PB-32189, Clearinghouse for Federal Scientific and
  Technical Information, U.S. Dept. of Commerce, Washington, D.C., 1942.

\bibitem[CF48]{courant-friedrichs}
{R.} Courant and {K.O.} Friedrichs, \emph{Supersonic flow and shock waves},
  Interscience Publishers, 1948.

\bibitem[Che20]{shuxing-chen-book-2020}
Shuxing Chen, \emph{Mathematical analysis of shock wave reflection}, Springer,
  2020.

\bibitem[EL08]{elling-liu-pmeyer}
{V.} Elling and Tai-Ping Liu, \emph{Supersonic flow onto a solid wedge}, Comm.
  Pure Appl. Math. \textbf{61} (2008), no.~10, 1347--1448.

\bibitem[Ell09]{elling-sonic-potf}
V.~Elling, \emph{Counterexamples to the sonic criterion}, Arch. Rat. Mech.
  Anal. \textbf{194} (2009), no.~3, 987--1010.

\bibitem[Ell21]{elling-idealpolar}
{V.} Elling, \emph{Shock polars for ideal and non-ideal gas}, J. Fluid Mech.
  \textbf{916} (2021), no.~A51.

\bibitem[Ell22a]{elling-isenpolar}
\bysame, \emph{Barotropic {Euler} shock polars}, Z. Angew. Math. Phys.
  \textbf{73} (2022).

\bibitem[Ell22b]{elling-polarpotf}
\bysame, \emph{Shock polars for non-polytropic compressible potential flow},
  Comm. Pure Appl. Anal. \textbf{21} (2022), no.~5, 1581--1594.

\bibitem[Ell22c]{elling-waalspolar}
\bysame, \emph{Van der {Waals} shock polars with multiple or supersonic
  critical points}, Phys. Fluids \textbf{34} (2022), 036110.

\bibitem[Fow81]{fowles-jfm1981}
{G.R.} Fowles, \emph{Stimulated and spontaneous emission of acoustic waves from
  shock fronts}, Phys. Fluids \textbf{24} (1981), 220--227.

\bibitem[FT68]{ferrari-tricomi}
{C.} Ferrari and {F.} Tricomi, \emph{Transonic aerodynamics}, Academic Press,
  1968.

\bibitem[HM98]{henderson-menikoff}
{L.F.} Henderson and {R.} Menikoff, \emph{Triple-shock entropy theorem and its
  consequences}, J. Fluid Mech. \textbf{366} (1998), 179--210.

\bibitem[HT02]{hunter-tesdall}
{J.} Hunter and {A.} Tesdall, \emph{Self-similar solutions for weak shock
  reflection}, {SIAM} J. Appl. Math. \textbf{63} (2002), no.~1, 42--61.

\bibitem[Jam78]{jameson-1978}
{A.} Jameson, \emph{Remarks on the calculation of transonic potential flow by a
  finite volume method}, Proceedings of the IMA Conference on Numerical Methods
  in Applied Fluid Mechanics ({B.} Hunt, ed.), Academic Press, January 1978,
  pp.~363--386.

\bibitem[K\'47]{karman-tsd-1947}
{T.}~{von} K\'{a}rm\'{a}n, \emph{The similarity law of transonic flow}, Journal
  of Mathematics and Physics \textbf{26} (1947), no.~1--4, 182--190.

\bibitem[KL13]{eunheuikim-chungminlee-2013}
Eun-Heui Kim and Chung-Min Lee, \emph{Numerical solutions to shock reflection
  and shock interaction problems for the self-similar transonic two-dimensional
  nonlinear wave systems}, Journal of Computational Science \textbf{4} (2013),
  no.~1, 36--45.

\bibitem[Lai19]{geng-lai-indiana2019}
{G.} Lai, \emph{Global solutions to a class of two-dimensional {Riemann}
  problems for the {Euler} equations with a general equation of state}, Indiana
  Univ. Math. J. \textbf{68} (2019), no.~5, 1409--1464.

\bibitem[LZY98]{li-zhang-yang}
Jiequan Li, Tong Zhang, and Shuli Yang, \emph{The two-dimensional {Riemann}
  problem in gas dynamics}, Addison Wesley Longman, 1998.

\bibitem[MC71]{murman-cole}
{E.M.} Murman and {J.D.} Cole, \emph{Calculation of plane steady transonic
  flows}, AIAA J. \textbf{9} (1971), no.~1, 114--121.

\bibitem[Pre38]{preiswerk-phdthesis1938}
{E.} Preiswerk, \emph{Anwendung gasdynamischer {Methoden} auf
  {Wasserstr\"omungen} mit freier {Oberfl\"ache}}, Ph.D. thesis, ETH
  {Z\"urich}, 1938, translated as NACA Technical Notes 934 and 935.

\bibitem[Tes86]{teshukov-polar}
{V.M.} Teshukov, \emph{On the shock polars in a gas with general equations of
  state}, J. Appl. Math. Mech. \textbf{50} (1986), no.~1, 71--75.

\bibitem[Wey49]{weyl-shock-waves}
H.~Weyl, \emph{Shock waves in arbitrary fluids}, Comm. Pure Appl. Math.
  \textbf{2} (1949), no.~2--3, 103--122.

\bibitem[Zha99]{yongqian-zhang-wedge-glimm}
Yongqian Zhang, \emph{Global existence of steady supersonic potential flow past
  a curved wedge with a piecewise smooth boundary}, {SIAM} J. Math. Anal.
  \textbf{31} (1999), 166--183.

\bibitem[Zhe01]{yuxi-zheng-book}
Yuxi Zheng, \emph{Systems of conservation laws}, {Birkh\"auser}, 2001.

\end{thebibliography}
\providecommand{\bysame}{\leavevmode\hbox to3em{\hrulefill}\thinspace}
\providecommand{\MR}{\relax\ifhmode\unskip\space\fi MR }
\providecommand{\MRhref}[2]{%
  \href{http://www.ams.org/mathscinet-getitem?mr=#1}{#2}
}
\providecommand{\href}[2]{#2}

\end{document}